\documentclass[a4paper,12pt, reqno]{amsart}
\usepackage[backref]{hyperref}
\usepackage{a4wide}
\usepackage{amsthm}

 \newtheorem{theorem}{Theorem}

\begin{document}

\title{Multiplicative inverses in short intervals}

\author{S. Baier}

\address{S. Baier, School of Mathematics, University of East Anglia, Norwich, NR4 7TJ, England}

\email{s.baier@uea.ac.uk}

\maketitle

\section{Introduction and results}
T.D. Browning and A. Haynes \cite{BrHa} considered the following problem. 
Let $p$ be a prime, $J$ be an integer and $I_1^{(j)},I_2^{(j)}$ with $1\le j\le J$ be finite sequences of subintervals of $(0,p)$. Under which conditions is there a $j$ such that  
\begin{equation} \label{congruence}
xy\equiv 1 \bmod{p}, \quad (x,y)\in \left(I_1^{(j)}\times I_2^{(j)}\right)\cap \mathbb{Z}^2
\end{equation}
has a solution? They proved the following.

\begin{theorem} \label{theorem1} Let $H,K>0$ and let $I_1^{(j)},I_2^{(j)} \subseteq (0,p)$ be subintervals, for $1\le j\le J$, such that
$$
\left|I_1^{(j)}\right|=H \quad \mbox{and} \quad \left|I_2^{(j)}\right|=K
$$
and
\begin{equation} \label{empty}
I_1^{(j)}\cap I_1^{(k)}=\emptyset \quad \mbox{for all} \quad j\not=k.
\end{equation}
Then there exists $j\in\{1,...,J\}$ for which \eqref{congruence} has a solution if
\begin{equation} \label{Jineq}
J\gg \frac{p^3\log^4 p}{H^2K^2}. 
\end{equation}
\end{theorem}

The proof of this theorem in \cite{BrHa} relies on the following new mean value theorem for short Kloosterman sums by Browning and Haynes.

\begin{theorem} \label{theorem2} If $I_1,...,I_J \subseteq (0,p)$ are disjoint subintervals, with $H/2\le |I_j|\le H$ for each $j$, then for any $l\in \left(\mathbb{Z}/p\mathbb{Z}\right)^{\ast}$, we have
$$
\sum\limits_{j=1}^J \left|\sum\limits_{n\in I_j} e\left(\frac{l\overline{n}}{p} \right)\right|^2 \le 2^{12}p\log^2 H.
$$
\end{theorem}

In this note, we prove Theorem \ref{theorem1}, in a slightly generalized and refined form, by a different method which doesn't use Theorem \ref{theorem2} but Poisson summation and Weil's estimate for Kloosterman sums. Moreover,  we improve this result under certain additional conditions on the spacing of the intervals $I_i^{(j)}$. 

We note that we could add the assumption 
\begin{equation} \label{Hass}
H\ge \log p
\end{equation}
to Theorem \ref{theorem1} without weakening the result because $J\ll p/H$ under the conditions of this theorem which contradicts \eqref{Jineq} if $H<\log p$ since $K\le p$. We want to assume \eqref{Hass} throughout the sequel.
Moreover, we want to assume without loss of generality that the intervals $I_1^{(j)}$ and $I_2^{(j)}$ are closed and centered at integers $N_j$ and $M_j$, respectively. 
We prove the following.

\begin{theorem} \label{theorem3} 
Assume that all conditions in Theorem \ref{theorem1} are satisfied, except possibly \eqref{empty}.  Assume the integers $M_j$ are  $X$-spaced modulo $p$, i.e.
$$
p\cdot \left|\left| \frac{M_j-M_k}{p} \right|\right| \ge X \quad \mbox{for all $j,k$ with } j\not=k,
$$
for some $X\ge 1$,
where $||z||$ denotes the distance of the real number $z$ to the nearest integer.
Then there exists $j\in \{1,...,J\}$ for which \eqref{congruence} has a solution if 
\begin{equation} \label{Jlowerbound}
J\gg \frac{p^3\log^{3+\varepsilon} p}{HK^2\min\{H,X\}}.
\end{equation}
\end{theorem}

Clearly, Theorem \ref{theorem3} implies Theorem \ref{theorem1}. If the intervals $I_1^{(j)}$ and $I_2^{(j)}$ are equispaced, respectively, we obtain an improvement of the result above, provided that the intervals $I_1^{(j)}$ are not spaced too far away. 

\begin{theorem} \label{theorem4}
Assume that all conditions in Theorem \ref{theorem1} are satisfied. Assume further that the integers $M_j$ and $N_j$ lie in arithmetic progression, respectively, i.e.
$$
M_j=M+jX \quad \mbox{and} \quad N_j=N+jY
$$
for certain integers $M,N,X,Y$ and all $j\in \{1,...,J\}$. Then there exists $j\in\{1,...,J\}$ for which \eqref{congruence} has a solution if
\begin{equation} \label{has}
X\ll \frac{HK}{p^{1/2}(\log p)^{1+\varepsilon}} \quad \mbox{and} \quad J\gg \frac{p^{3/2}\log^{2+\varepsilon} p}{HK}.
\end{equation}
\end{theorem}

\section{Basic approach}
Set
$$
w(t):=\exp\left(-\pi t^2\right).
$$
This has Fourier transform
\begin{equation} \label{Fourier}
\hat{w}(t)=w(t)=\exp\left(-\pi t^2\right).
\end{equation}
Set 
\begin{equation} \label{6}
x:=\frac{H}{(\log p)^{1/2+\varepsilon}}, \quad y:=\frac{K}{(\log p)^{1/2+\varepsilon}}
\end{equation}
and
$$
T:=\sum\limits_{j=1}^J \mathop{\sum\limits_{|m-M_j|\le H/2} \sum\limits_{|n-N_j|\le K/2}}_{m\equiv \overline{n} \bmod{p}} w\left(\frac{m-M_j}{x}\right)  w\left(\frac{n-N_j}{y}\right),
$$
where here and in the sequel, we assume that $(n,p)=1$, and $\overline{n}$ denotes the multiplicative inverse of $n$ modulo $p$.
Then clearly, \eqref{congruence} has a solution if  $T>0$.
Now the general strategy is to extend the sums over $m$ and $n$ to all integers and use Poisson summation and Weil's estimate for Kloosterman sums.

We write
\begin{equation} \label{write}
T=S-S_1-S_2,
\end{equation}
where 
$$
S:=\sum\limits_{j=1}^J \mathop{\sum\limits_{m} \sum\limits_{n}}_{m\equiv \overline{n} \bmod{p}} w\left(\frac{m-M_j}{x}\right)  w\left(\frac{n-N_j}{y}\right),
$$
$$
S_1:=\sum\limits_{j=1}^J \mathop{\sum\limits_{|m-M_j|> H/2} \sum\limits_{n}}_{m\equiv \overline{n} \bmod{p}} w\left(\frac{m-M_j}{x}\right)  w\left(\frac{n-N_j}{y}\right)
$$
and 
$$
S_2:=\sum\limits_{j=1}^J \mathop{\sum\limits_{|m-M_j|\le H/2} \sum\limits_{|n-N_j|>K/2}}_{m\equiv \overline{n} \bmod{p}} w\left(\frac{m-M_j}{x}\right)  w\left(\frac{n-N_j}{y}\right).
$$
From $J,H,K\le p$ ($J\le p$ following from $X\ge 1$), \eqref{Fourier} and \eqref{6}, it is evident that  $S_1$ and $S_2$ are negligible, i.e.
\begin{equation} \label{negligible}
S_1,S_2\ll_A p^{-A} \quad \mbox{for any } A>0.
\end{equation}
In the following section, we estimate the sum $S$.

\section{Application of Poisson summation}
We have 
\begin{equation*} \label{trans}
\begin{split}
S= & \sum\limits_{j=1}^J \mathop{\sum\limits_{m}  \sum\limits_{n}}_{m\equiv \overline{n} \bmod{p}} 
w\left(\frac{m-M_j}{x}\right)w\left(\frac{n-N_j}{y}\right)
\\
= & \frac{1}{p} \sum\limits_{k=-(p-1)/2}^{(p-1)/2} \sum\limits_{j=1}^J \sum\limits_{m} \sum\limits_{n} w\left(\frac{m-M_j}{x}\right) w\left(\frac{n-N_j}{y}\right)
e\left(k\cdot \frac{m-\overline{n}}{p}\right) \\
= & \frac{1}{p} J \sum\limits_{m} w\left(\frac{m}{x}\right)\sum\limits_{n} w\left(\frac{n}{y}\right)+\\ &  \frac{1}{p} \sum\limits_{\substack{k=-(p-1)/2\\ k\not=0}}^{(p-1)/2} \sum\limits_{j=1}^J \sum\limits_{m} w\left(\frac{m-M_j}{x}\right)e\left(k\cdot \frac{m}{p}\right)\sum\limits_{n} w\left(\frac{n-N_j}{y}\right)e\left(-k\cdot \frac{\overline{n}}{p}\right).
\end{split}
\end{equation*}
Using Poisson summation, the terms in the last line can be transformed as follow. First,
\begin{equation*} \label{poisson1}
\frac{1}{p} J \sum\limits_{m} w\left(\frac{m}{x}\right)\sum\limits_{n} w\left(\frac{n}{y}\right)= \frac{Jxy}{p} \cdot \hat{w}(0)^2. 
\end{equation*}
Second,
\begin{equation*} \label{poisson2}
\begin{split}
& \sum\limits_{m} w\left(\frac{m-M_j}{x}\right)e\left(k\cdot \frac{m}{p}\right) = \sum\limits_{m} w\left(\frac{m}{x}\right) e\left(k\cdot \frac{m+M_j}{p}\right)\\ =& 
x\cdot e\left(k\cdot \frac{M_j}{p}\right) \cdot \sum\limits_{m\equiv k \bmod{p}} \hat{w}\left(\frac{mx}{p}\right)=x\cdot e\left(k\cdot \frac{M_j}{p}\right) \cdot F_k(x), \quad \mbox{say.}
\end{split}
\end{equation*} 
Third,
\begin{equation*} \label{poitsson3}
\begin{split}
& \sum\limits_{n} w\left(\frac{n-N_j}{y}\right)e\left(-k\cdot \frac{\overline{n}}{p}\right) 
= \sum\limits_{c\bmod{p}} e\left(-k\cdot \frac{\overline{c}}{p}\right) \sum\limits_{n \equiv c+N_j\bmod{p}} w\left(\frac{n}{y}\right)\\
=& \frac{y}{p}\cdot  \sum\limits_{c\bmod{p}} e\left(-k\cdot \frac{\overline{c}}{p}\right) \sum\limits_{l} \hat{w}\left(\frac{ly}{p}\right)e\left(l\cdot \frac{c+N_j}{p}\right)\\
=& \frac{y}{p} \cdot \sum\limits_{l} \hat{w}\left(\frac{ly}{p}\right) e\left(l\cdot \frac{N_j}{p}\right) S(l,-k;p),
\end{split}
\end{equation*}
where 
$$
S(l,-k;p)=\sum\limits_{c=1}^{p-1} e\left(\frac{lc-k\overline{c}}{p}\right)
$$
is the Kloosterman sum.
Putting the above together, we get
\begin{equation}\label{above}
\begin{split}
S= & \frac{Jxy}{p} \cdot \hat{w}(0)^2+\\ & \frac{xy}{p^2} \cdot   \sum\limits_{\substack{k=-(p-1)/2\\ k\not=0}}^{(p-1)/2} \sum\limits_{l} S(l,-k;p) F_k(x) \hat{w}\left(\frac{ly}{p}\right) \sum\limits_{j=1}^J e\left(\frac{kM_j+lN_j}{p}\right).
\end{split}
\end{equation}

We note that if $-(p-1)/2\le k\le (p-1)/2$, then
\begin{equation} \label{Fkx}
\begin{split}
F_k(x)= & \sum\limits_{m\equiv k \bmod{p}} \hat{w}\left(\frac{mx}{p}\right) =  \sum\limits_{r\in \mathbb{Z}} \hat{w}\left(\left(\frac{k}{p}+r\right)x\right)\\
\ll & \exp\left(-\frac{|k|}{p}\cdot x\right) \cdot \sum\limits_{r=0}^{\infty} \exp(-rx) \ll  \exp\left(-\frac{|k|}{p}\cdot x\right)
\end{split}
\end{equation}
since $x\ge 1$ if $\varepsilon\le 1/2$ by \eqref{Hass} and \eqref{6}.

\section{Proof of Theorem 3}
By Weil's bound for Kloosterman sums, we have 
\begin{equation} \label{Weil}
S(l,-k;p)\ll p^{1/2} \quad \mbox{if } k\not\equiv 0 \bmod{p}.
\end{equation}
Using the Cauchy-Schwarz inequality and \eqref{Fkx}, it follows that
\begin{equation} \label{aftercauchy}
\begin{split}
& \sum\limits_{\substack{k=-(p-1)/2\\ k\not=0}}^{(p-1)/2} \sum\limits_{l} S(l,-k;p) F_k(x) \hat{w}\left(\frac{ly}{p}\right) \sum\limits_{j=1}^J e\left(\frac{kM_j+lN_j}{p}\right)\\ \ll &
p^{1/2}\left(\sum\limits_{\substack{k=-(p-1)/2\\ k\not=0}}^{(p-1)/2} \sum\limits_{l} F_k(x) \hat{w}\left(\frac{ly}{p}\right)\right)^{1/2}  \\ &
\left(\sum\limits_{k=-(p-1)/2}^{(p-1)/2} \sum\limits_{l}  F_k(x) \hat{w}\left(\frac{ly}{p}\right)\left|\sum\limits_{j=1}^J e\left(\frac{kM_j+lN_j}{p}\right)\right|^2\right)^{1/2}\\
\ll & \frac{p^{3/2}}{(xy)^{1/2}} \left(\sum\limits_{k=-(p-1)/2}^{(p-1)/2} \sum\limits_{l}  F_k(x) \hat{w}\left(\frac{ly}{p}\right)\left|\sum\limits_{j=1}^J e\left(\frac{kM_j+lN_j}{p}\right)\right|^2\right)^{1/2}.
\end{split}
\end{equation}
Expanding the square, we get
\begin{equation*}
\begin{split}
& \sum\limits_{k=-(p-1)/2}^{(p-1)/2} \sum\limits_{l}  F_k(x) \hat{w}\left(\frac{ly}{p}\right)\left|\sum\limits_{j=1}^J e\left(\frac{kM_j+lN_j}{p}\right)\right|^2\\ =& \sum\limits_{j_1,j_2=1}^J \left(\sum\limits_{k=-(p-1)/2}^{(p-1)/2} F_k(x)e\left(k\cdot \frac{M_{j_1}-M_{j_2}}{p}\right)\right)
\left(\sum\limits_{l} \hat{w}\left(\frac{ly}{p}\right)e\left(l\cdot \frac{N_{j_1}-N_{j_2}}{p}\right)\right).\end{split}
\end{equation*}
Using \eqref{Fkx}, we have 
$$
\sum\limits_{k=-(p-1)/2}^{(p-1)/2} F_k(x)e\left(k\cdot \frac{M_{j_1}-M_{j_2}}{p}\right) \ll \sum\limits_{k=-(p-1)/2}^{(p-1)/2} F_k(x) \ll \frac{p}{x}.
$$
Similarly,
$$
\sum\limits_{l} \hat{w}\left(\frac{ly}{p}\right)e\left(l\cdot \frac{N_{j_1}-N_{j_2}}{p}\right) \ll \sum\limits_{l} \hat{w}\left(\frac{ly}{p}\right) \ll \frac{p}{y}.
$$
Moreover, removing the weight function $F_k(x)$ using partial summation and using the familiar estimate for geometric sums, we have
$$
\sum\limits_{k=-(p-1)/2}^{(p-1)/2} F_k(x)e\left(k\cdot \frac{M_{j_1}-M_{j_2}}{p}\right) \ll \left|\left| \frac{M_{j_1}-M_{j_2}}{p}\right|\right|^{-1}.
$$
It follows that  
\begin{equation} \label{ale}
\begin{split}
& \sum\limits_{k=-(p-1)/2}^{(p-1)/2} \sum\limits_{l} F_k(x) \hat{w}\left(\frac{ly}{p}\right)\left|\sum\limits_{j=1}^J e\left(\frac{kM_j+lN_j}{p}\right)\right|^2\\
\ll & \frac{p}{y}\sum\limits_{j_1,j_2=1}^J\min\left\{\left|\left| \frac{M_{j_1}-M_{j_2}}{p}\right|\right|^{-1}, \frac{p}{x}\right\} .
\end{split}
\end{equation}
Using the fact that the $M_j$'s are $X$-spaced modulo $p$, we obtain 
\begin{equation} \label{last}
\begin{split}
\ll & \sum\limits_{j_1,j_2=1}^J\min\left\{\left|\left| \frac{M_{j_1}-M_{j_2}}{p}\right|\right|^{-1}, \frac{p}{x}\right\} \\
\ll & J \sum\limits_{j=0}^{J-1} \min\left\{\frac{p}{jX}, \frac{p}{x}\right\}\\
\ll & Jp\left(\frac{1}{x}+\frac{\log 2J}{X}\right). 
\end{split}
\end{equation}

Combining \eqref{write}, \eqref{negligible}, \eqref{above}, \eqref{aftercauchy}, \eqref{ale} and \eqref{last}, we arrive at
\begin{equation}
\begin{split}
T = \frac{Jxy}{p} \cdot \hat{w}(0)^2+O\left((\log 2J)^{1/2}(Jxp)^{1/2}
\left(\frac{1}{x}+\frac{1}{X}\right)^{1/2}\right).
\end{split}
\end{equation}
For the right-hand side to be greater $0$ (i.e. error term $<$ main term), it suffices that
$$
J\gg \frac{p^3\log p}{xy^2}\left(\frac{1}{x}+\frac{1}{X}\right),
$$
which holds if \eqref{Jlowerbound} is satisfied. This implies Theorem \ref{theorem3}.

\section{Proof of Theorem 4}
Under the conditions of Theorem 4, we have 
\begin{equation} \label{geom}
\begin{split}
\sum\limits_{j=1}^J e\left(\frac{kM_j+lN_j}{p}\right)=e\left(\frac{kM+lN}{p}\right) \sum\limits_{j=1}^J e\left(\frac{j(kX+lY)}{p}\right)\\ \ll 
\min\left\{\left|\left|\frac{kX+lY}{p}\right|\right|^{-1},J\right\}.
\end{split}
\end{equation}
Using \eqref{Weil} and \eqref{geom}, we obtain
\begin{equation} \label{follows}
\begin{split}
 & \sum\limits_{k=-(p-1)/2}^{(p-1)/2} \sum\limits_{l} S(l,-k;p) F_k(x) \hat{w}\left(\frac{ly}{p}\right) \sum\limits_{j=1}^J e\left(\frac{kM_j+lN_j}{p}\right)\\
 \ll & p^{1/2}  \sum\limits_{l}  \hat{w}\left(\frac{ly}{p}\right) \sum\limits_{k=-(p-1)/2}^{(p-1)/2}F_k(x) \min\left\{\left|\left|\frac{kX+lY}{p}\right|\right|^{-1},J\right\}.
\end{split}
\end{equation}
We estimate the inner-most sum over $k$ by
\begin{equation*} 
\begin{split}
\sum\limits_{k=-(p-1)/2}^{(p-1)/2}F_k(x)\min\left\{\left|\left|\frac{kX+lY}{p}\right|\right|^{-1},J\right\} \ll & \frac{p/x}{p/X}\cdot \left(J+\frac{p}{X}\cdot\log p\right)\\
= & \frac{X}{x}\cdot J+ \frac{p}{x}\cdot \log p, 
\end{split}
\end{equation*}
where we use \eqref{Fkx} and $X\ge x$ (which follows from \eqref{empty}). Hence, we get
\begin{equation} \label{hence}
\begin{split}
& p^{1/2}  \sum\limits_{l}  \hat{w}\left(\frac{ly}{p}\right) \sum\limits_{k=-(p-1)/2}^{(p-1)/2}F_k(x) \min\left\{\left|\left|\frac{kX+lY}{p}\right|\right|^{-1},J\right\}\\
\ll & \frac{p^{3/2}}{y} \cdot \left(\frac{X}{x}\cdot J+ \frac{p}{x}\cdot \log p \right).
\end{split}
\end{equation}
Combining \eqref{write}, \eqref{negligible}, \eqref{above}, \eqref{follows} and \eqref{hence}, we get
\begin{equation}
T=\frac{Jxy}{p} \cdot \hat{w}(0)^2+O\left(\frac{X}{p^{1/2}}\cdot J+p^{1/2}\log p \right).
\end{equation}
For the right-hand side to be greater 0, it suffices that
$$
X\ll \frac{xy}{p^{1/2}} \quad \mbox{and} \quad J\gg \frac{p^{3/2}\log p}{xy},
$$
which holds if \eqref{has} is satisfied. This implies Theorem \ref{theorem4}.

\end{document}